\documentclass[graybox]{svmult}

\usepackage{mathptmx}       % selects Times Roman as basic font; note that this font does
\usepackage{helvet}         % selects Helvetica as sans-serif font
\usepackage{courier}        % selects Courier as typewriter font
\usepackage{graphicx}              % standard LaTeX graphics tool when including figure files

\usepackage{siunitx}
\usepackage{amsbsy}
\usepackage{amsmath,bm}
\usepackage{amssymb}
\usepackage{amsfonts}
\usepackage{algorithmic}
\usepackage{algorithm} % !!do NOT use the package algorithm2e as it will clash with the algorithm package!!

\usepackage{multirow}
\usepackage[bottom]{footmisc}% places footnotes at page bottom
\usepackage{cite}
\usepackage{url}

\begin{document}

\title*{Towards topology optimization of a hybrid - excited machine using recursive material interpolation}
\titlerunning{Multimaterial topology optimization using recursive material interpolation}
\author{Théodore Cherrière
}
\institute{Théodore Cherrière \at Johann Radon Institute for Computational and Applied Mathematics (RICAM), Austrian Academy of Science, A-4040 Linz, Altenbergerstraße 69,
\email{theodore.cherriere@ricam.oeaw.ac.at}.\\
This work is associated with FWF-funded CRC F90 (CREATOR – Computational Electric Machine Laboratory), project 08.
}

\maketitle

\abstract*{Hybrid-excited electrical machines aim to combine the advantages of permanent magnet machines (high efficiency and torque density) with those of separately excited machines (ease of flux-weakening at high speed). These machines are of interest to electric vehicles, and only parametric approaches are available in the literature for their optimization. This work proposes a more general topology optimization methodology by extending the formalism of density methods. The difficulty lies in integrating the numerous natures of materials (conductors, permanent magnets, ferromagnetic material, air...) without strongly deconvexifying the optimization problem, which leads to non-physical results with unsatisfactory performance. To address this issue, a recursive material interpolation is introduced. The hybrid-excited rotors optimized by this approach are compared with those of existing techniques, demonstrating a clear superiority of the recursive interpolation.}

\abstract{Hybrid-excited electrical machines aim to combine the advantages of permanent magnet machines (high efficiency and torque density) with those of separately excited machines (ease of flux-weakening at high speed). These machines are of interest to electric vehicles, and only parametric approaches are available in the literature for their optimization. This work proposes a more general topology optimization methodology by extending the formalism of density methods. The difficulty lies in integrating the numerous natures of materials (conductors, permanent magnets, ferromagnetic material, air...) without strongly deconvexifying the optimization problem, which leads to non-physical results with unsatisfactory performance. To address this issue, a recursive material interpolation is introduced. The hybrid-excited rotors optimized by this approach are compared with those of existing techniques, demonstrating a clear superiority of the recursive interpolation.}

\section{Introduction}
Topology optimization is a promising tool for designing innovative electromagnetic actuators, overcoming state-of-the-art structures. Unlike conventional optimization methodologies relying on a geometric parametrization of an initial structure arbitrarily chosen by the designer from experience \cite{Cherriere:Cisse2021}, topology optimization is nonparametric and thus aims to extend the space of possible results. Initially developed in mechanical engineering \cite{Cherriere:Bendsoe1989}, topology optimization was progressively applied to electromagnetic actuators from the 1990s \cite{Cherriere:Dyck1996}.  Topology optimization is now mature for single-material problems and is implemented in many commercial software.

\begin{figure}
	\centering
	\includegraphics[height=2.1cm]{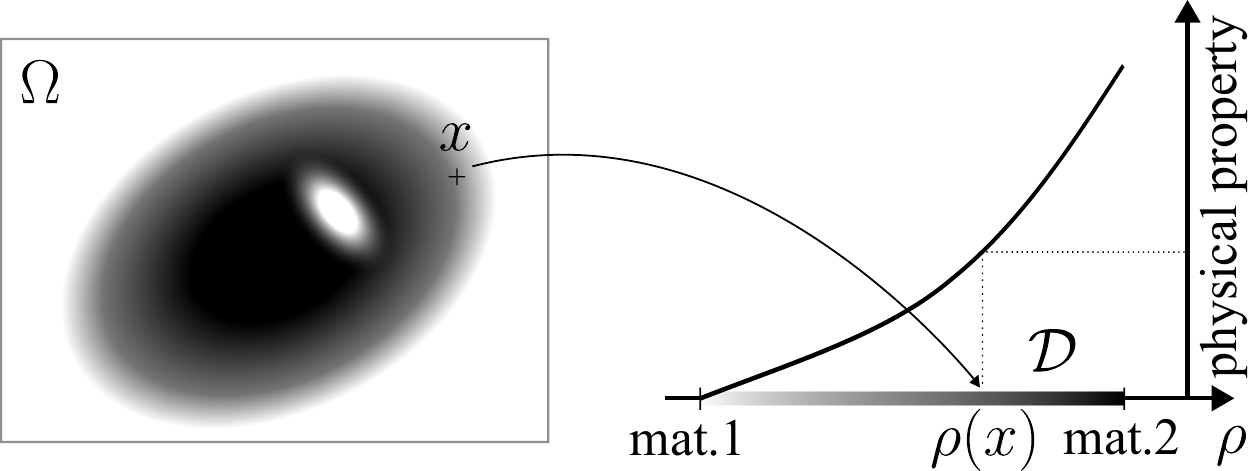}\qquad
	\includegraphics[height=2.1cm]{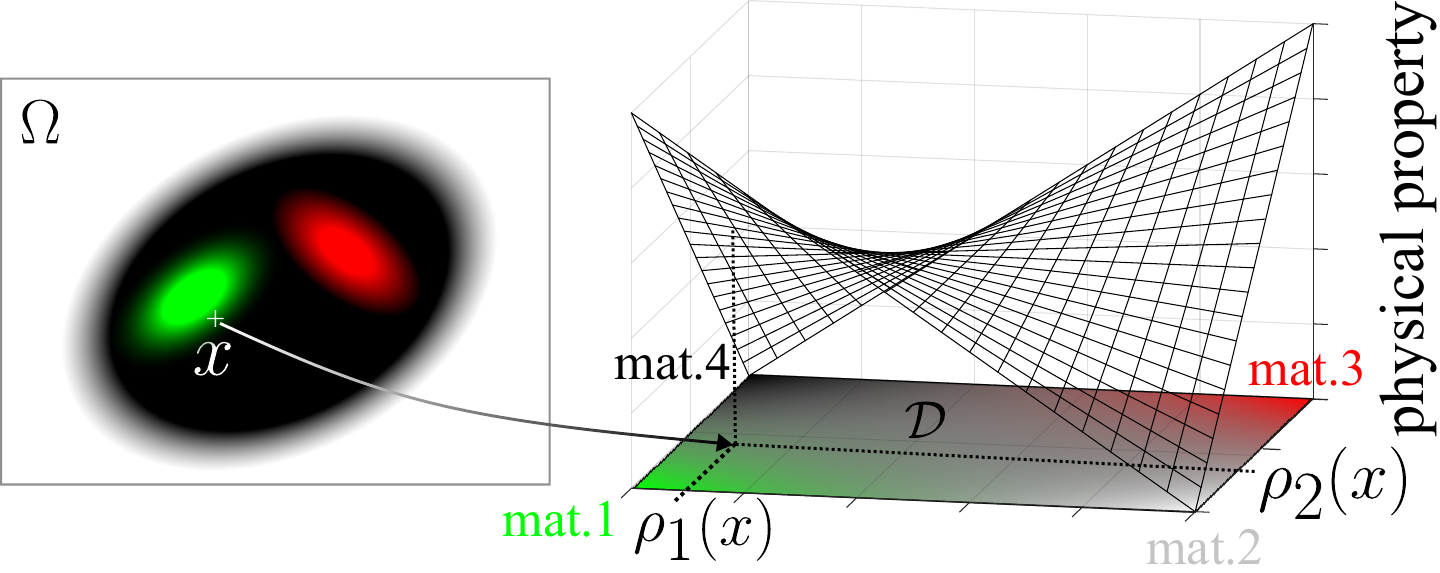}
	\caption{Interpolation of single/bi-material (left) and multi-material (right) with density methods.} \label{fig:Cherriere_principleDensite}
\end{figure}

Among the various topology optimization techniques (\textit{cf} \cite{Cherriere:Lucchini2022} for a detailed state of the art), density methods are the most widely used. They rely on a so-called "density" field $\rho: \Omega \rightarrow \mathcal{D} $ over the design region $x\in \Omega$, from which the physical properties at each point $x\in \Omega$ are defined by interpolating the candidate materials' properties. Each vertex of the interpolation domain $\mathcal{D}$ is associated with a candidate material, so to consider more than two materials, the dimension of $\mathcal{D}$ must be increased. The literature generally considers hypercubic domains \cite{Cherriere:Bruyneel2011} so that $\mathcal{D} = [0,1]^n$, with $n\in \mathbb{N}^*$, as illustrated in Figure~\ref{fig:Cherriere_principleDensite} for $n=1$ and $2$. Note that despite its name, there is generally no physical interpretation of intermediate values of $\rho$.

We consider an objective function $f$ to be maximized, which generally depends on a physical state $a\in H^1_0(\Omega)$ that implicitly depends on the design field $\rho$ by solving the physical equations. The optimization problem is then written as
\begin{equation}
            \text{find} \quad \quad  \rho_\text{opt} = \arg \displaystyle   \min_{\rho : \Omega \rightarrow \mathcal{D}} f(a(\rho)). \label{eq:Cherriere_PbOptim}
\end{equation}

Since the sensitivity $\mathrm{d}_\rho f$ can be computed efficiently by the adjoint method \cite{Cherriere:Cea1986}, we can solve \eqref{eq:Cherriere_PbOptim} by a projected gradient descent \cite{Cherriere:Calamai1987}. Unfortunately, hypercubic domains are not suitable for all optimization problems since the resulting interpolation of physical properties is often non-monotonic, as shown in Figure \ref{fig:Cherriere_principleDensite}, which strongly dexonvexifies the problem and leads to uninteresting local optima. One solution is to extend $\mathcal{D}$ to the set of convex polytopes \cite{Cherriere:Cherriere2022}. These polytopes allow materials of the same nature (typically, conductors or magnets with different phases or orientations) to lie on the same plane and the other materials to be placed on an orthogonal axis.
However, this solution is limited in practice to low dimensions ($\dim \mathcal{D} \leq 3$) because of the complexity of the projection, shape function computation, and the difficulty of visualization to select appropriate higher-dimensional domains. This dimension limitation prevents more materials from being considered, as in the case of hybrid-excited electric machines that aim to combine the flexibility of wound-field excitation with the efficiency of permanent magnets (PM) for traction applications \cite{Cherriere:Zhu2019}.

In this context, this work presents a generalization of the interpolation framework required for multi-material density-based topology optimization. It relies on recursive interpolations supported on subdomains up to three dimensions. The implementation is freely available \cite{Cherriere:Cherriere2024Code} and is briefly described in Section~\ref{sec:Cherriere:Interpolation}. This formalism is then applied to the topology optimization of a hybrid-excited machine rotor and compared with existing techniques in the literature in Section~\ref{sec:Cherriere:TopOpt}.

\section{Interpolation framework}\label{sec:Cherriere:Interpolation}

Let us first consider $n_m$ candidate materials distributed over a single polytope $\mathcal{D}$. Each material $i$ is associated with the $i$-vertex of $\mathcal{D}$. The interpolation of a material property $\kappa$, possibly dependent on the physical state $a$ in the case of non-linear materials, is denoted $\tilde{\kappa}(a,\rho)$. This interpolation is constructed using generalized barycentric coordinates of $\mathcal{D}$ \cite{Cherriere:Cherriere2022}. We denote $\omega_i(\rho)$ the barycentric coordinate associated with vertex $i$, which is 1 when $\rho$ corresponds to the coordinates of vertex $i$, and 0 when $\rho$ corresponds to the coordinates of the other vertices. These barycentric coordinates partition the unit, and their computation is automatic \cite{Cherriere:Floater2014}. Therefore, a generic interpolation is defined by:
\begin{equation}
    \tilde{\kappa}(\rho,a) = \sum_{i=1}^{n_m} \omega_i(\rho) \kappa_i(a).
    \label{eq:Cherriere:interp}
\end{equation}

\begin{figure}
	\centering
	\includegraphics[height=4cm]{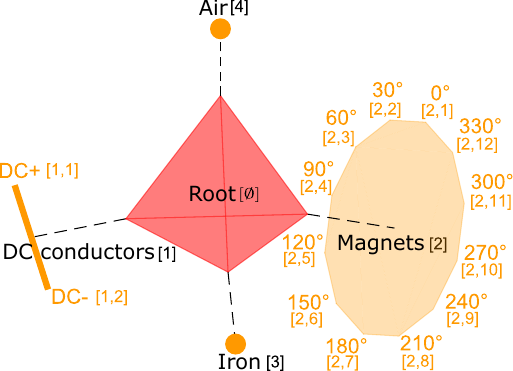}\qquad
	\includegraphics[height=4cm]{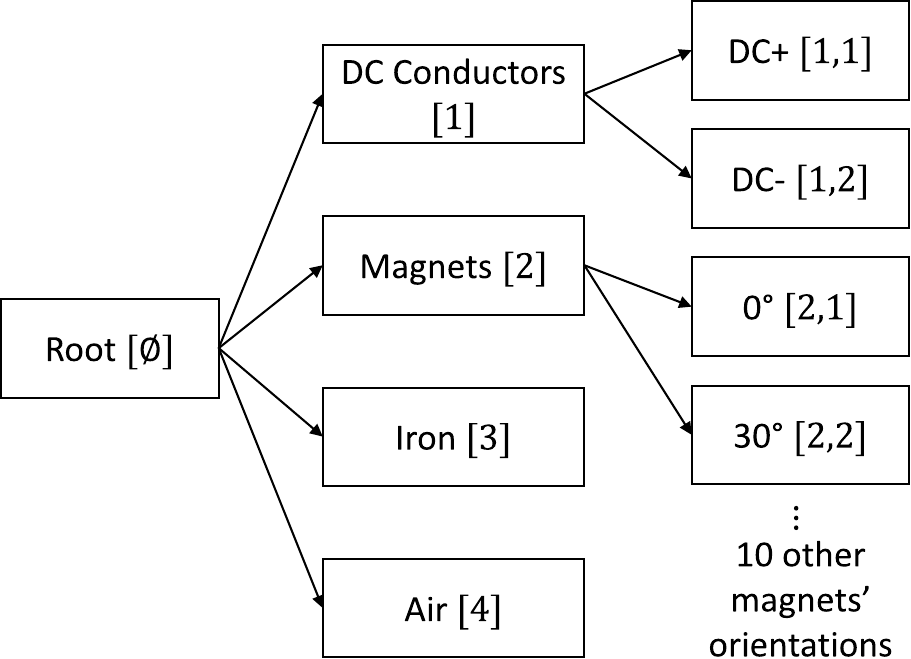}
	\caption{Hierarchical $\mathcal{D}$: 3D visualization (left) and associated rooted tree with Neveu's labels (right)} \label{fig:Cherriere:multiDomain}
\end{figure}

In \eqref{eq:Cherriere:interp}, another interpolation $\tilde{\kappa}_i$ on a new sub-domain can replace any material property $\kappa_i$. Figure~\ref{fig:Cherriere:multiDomain} (left) gives a convenient representation of a possible hierarchical interpolation domain. Projection onto such a structure is simply reduced to projection onto each low-dimensional sub-domain, which is easy to handle. The optimization variables are then dissociated and defined on each sub-domain.

From a computational point of view, this structure is a rooted tree denoted as $\mathcal{T}$, the branches being products by shape functions, the nodes being sums, and the leaves being material properties. Using Neveu's notation \cite{Cherriere:Neveu1986}, the set of children associated to the node $l$ is $\mathcal{C}(l)=\{[l,n]~|~n \in \mathbb{N}^*,[l,n]\in \mathcal{T}\}$, with $[\alpha,\beta]$ the concatenation of the lists of indices $\alpha$ and $\beta$. An example is shown in Figure~\ref{fig:Cherriere:multiDomain} (right).

Each subdomain $\mathcal{D}_l$ is associated with a specific optimization variable $\rho_l$. The evaluation of the interpolated physical property $\tilde{\kappa}_{\emptyset}$ %$\footnote{It is necessary to calculate the derivatives with respect to $\rho$ to compute the sensitivities of the objective function within the framework of a gradient descent, and the derivatives with respect to $a$ to solve the nonlinear physical problem with Newton's method (\textit{e.g.} by finite elements).}
can be performed recursively by applying the chain rule with $l=\emptyset$:

\begin{equation}
        \tilde{\kappa}_l(a,\rho) = \left \{
\begin{array}{ll}
 \kappa_l(a) & \quad \text{if~} \mathcal{C}(l) = \emptyset\\
\sum_{i\in \mathcal{C}(l)} \omega_i(\rho_l) \tilde{\kappa}_i(a,\rho) & \quad \text{else.}
\end{array}
\right. 
\label{eq:Cherriere:recursiveInterp}
\end{equation}

It is necessary to compute the derivatives with respect to $a$ to solve the nonlinear physical problem using Newton's method and the derivative with respect to $\rho$ to compute the sensitivities of  $f$ within the gradient descent framework. The evaluation of $\partial_a \tilde{\kappa}$ is straightforward by replacing $\tilde{\kappa}$ and $\kappa$ by $\partial_a \tilde{\kappa}$ and $\partial_a \kappa$ in \eqref{eq:Cherriere:recursiveInterp}, respectively. The computation of $\partial_\rho \tilde{\kappa}$ is a little bit more involved since the optimization variables are dissociated in each subdomain $\mathcal{D}_l$. Each partial derivative $\partial_{\rho_l} \tilde{\kappa}$ can therefore be evaluated recursively for all $l\in \mathcal{T}$ by calling \texttt{drho\_k}$(\emptyset,1)$ defined in Algorithm \ref{algo:Cherriere:RhoDerivative}.

\begin{algorithm}
\caption{\texttt{drho\_k}($l,k$) : Evaluate $\partial_\rho \tilde{\kappa}$}
\label{algo:Cherriere:RhoDerivative}
\begin{algorithmic}
%\REQUIRE $l$ optional Neveu's index, default: $\emptyset$ (corresponding to the root node)
%\REQUIRE $k\in\mathbb{R}$ optional constant, default: 1
%\REQUIRE $y$ optional hash table, default: hash table with no key
%\REQUIRE $\mathcal{T}$ (rooted tree), $\rho$ (hash-table of the optimization variables, callable with a Neveu's index of $\mathcal{T}$), $a$ (physical state), as well as $\kappa$ (hash-table of the candidate material properties), are assumed to be known global variables.
\STATE Compute $\mathcal{C}(l)$
\STATE $\displaystyle \partial_{\rho_l} \tilde{\kappa} \leftarrow k \sum_{i\in \mathcal{C}(l)}\tilde{\kappa}_l(a,\rho) \frac{\mathrm{d}\omega_i}{\mathrm{d}\rho_l}$ \quad \COMMENT{$\tilde{\kappa}_l(a,\rho)$ is computed with \eqref{eq:Cherriere:recursiveInterp}}
\FOR{$i\in \mathcal{C}(l)$}
\STATE \texttt{drho\_k}($i,k \omega_i(\rho_l)$)
\ENDFOR
\end{algorithmic}
\end{algorithm}

A general open-source Matlab$^\text{\textregistered}$ implementation is provided in \cite{Cherriere:Cherriere2024Code}, including elements for defining interpolation domains and functions to be interpolated, as well as evaluation, derivation, and projection routines. An example of application to electrical engineering is given in the next section. This implementation also allows for the consideration of penalty functions, common in topology optimization to eliminate intermediate materials \cite{Cherriere:Bendsoe1999}, but not used in this paper.

\section{Topology optimization of a hybrid-excited rotor} \label{sec:Cherriere:TopOpt}

\subsection{Problem definition}

The topology optimization of rotors for hybrid-excited machines \cite{Cherriere:Zhu2019} is a typical problem where several materials' natures coexist. These machines combine ferromagnetic material, air, PMs, and adjustable direct-current (DC) electrical conductors. The roles of the DC supply are (i) to reinforce the flux created by the magnets at low speeds to generate greater torque and (ii) to weaken the magnets' flux at high speeds to reduce the back-electromotive force.

\begin{figure}
	\centering
	\begin{tabular}{ccc}
		Hexadecagon (2D) & Diamond (3D) & Recursive domain (6D) \\ \hline \\
		\includegraphics[height=2.7cm]{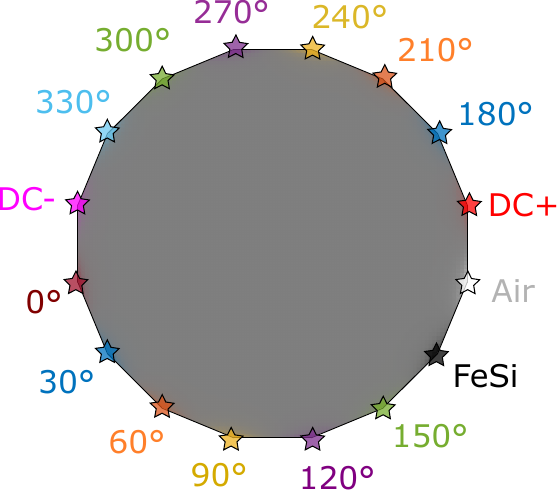} 
		& \includegraphics[height=2.7cm]{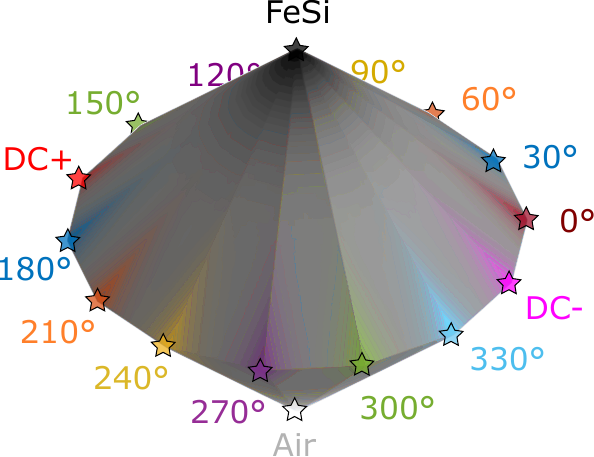}
		& \includegraphics[height=2.7cm]{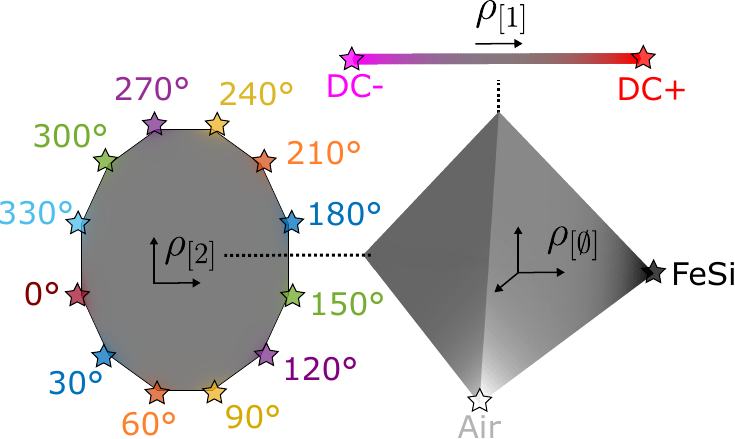} \\
	\end{tabular}
	\caption{3D representations of the interpolation domains with the materials' color scale.} \label{fig:Cherriere:domainApplication}
\end{figure}

We consider a rotor pole with 16 different candidate "materials": one ideal magnet with a remanent field of \SI{1}{T} and 12 discrete orientations, two electrical conductors with opposite current density equal to $\pm\SI{10}{A/mm^2}$, iron-silicon steel with standard anhysteretic nonlinear behavior \cite{Cherriere:Cherriere2023_M} and air. Several interpolation domains $\mathcal{D}$ are compared and plotted in Figure~\ref{fig:Cherriere:domainApplication}, one being defined recursively with the formalism described in Section~\ref{sec:Cherriere:Interpolation}.
The optimization problem reads:
\begin{equation}
         \text{find} \quad \displaystyle \rho^* = \min_{\rho \in \mathcal{D}} - \left( \frac{\gamma +1}{2}\phi^+(a(\rho)) -  \frac{\gamma-1}{2}  \phi^-(a(\rho)) \right),
         \label{eq:Cherriere_pbHybridFlux}
\end{equation}

\noindent with $\gamma\in[-1,1]$ an arbitrary coefficient, $\phi^+$ and $\phi^-$ the outward magnetic fluxes produced by the rotor pole when the supply current is positive and negative, respectively. The physical state $a$ is here the component along $z$ of the magnetic vector potential, which verifies the classical equations of 2D magnetostatics \cite{Cherriere:Cherriere2023_M}. The physical quantities interpolated by \eqref{eq:Cherriere:interp} are magnetic polarization and current density. The optimization region is a rotor pole of the same size as that of the BMWi3 motor, discretized on a mesh of 9663 elements.
The stator influence is neglected and is replaced by air. The optimization algorithm is a regularized projected gradient descent with density filtering that is detailed in \cite{Cherriere:Cherriere2024Compumag}. It stops after 500 iterations or if the design variables stagnate when the norm of the relative difference between two consecutive iterations is less than $10^{-4}$.

\subsection{Optimization results}

 Optimizations are performed for a regular distribution of $\gamma\in[-1,1]$ and homogeneous initializations in the different interpolation domains shown in Figure~\ref{fig:Cherriere:domainApplication}. 
 
According to \eqref{eq:Cherriere_pbHybridFlux}, the outwards flux should be maximized when $\gamma= 1$, and the inwards flux should be maximized when $\gamma= -1$. Then, the associated optimized rotors contain only PMs arranged in a well-known Halbach structure \cite{Cherriere:Zhu2001}, because the magnets generate higher flux than electrical conductors with reasonable current.

\begin{figure}
	\centering
	\begin{tabular}{ccc}
		$\gamma = -1$ & $\gamma = 0$ & $\gamma = 1$ \\ \hline \\
		\includegraphics[width=0.24\textwidth]{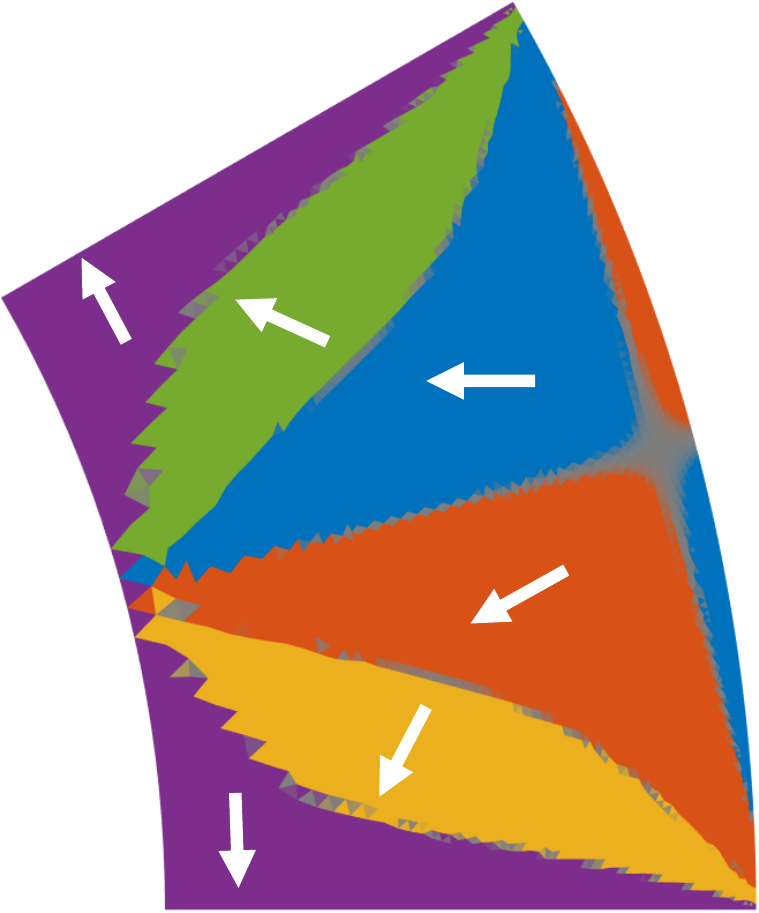} \qquad 
		& \qquad \includegraphics[width=0.24\textwidth]{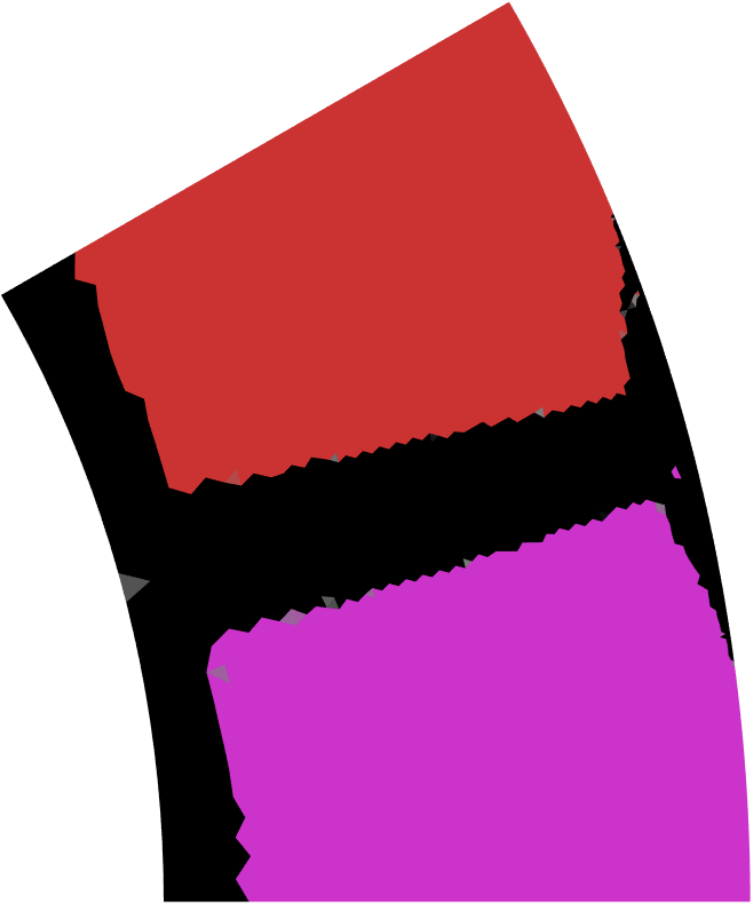}  \qquad
		& \qquad \includegraphics[width=0.24\textwidth]{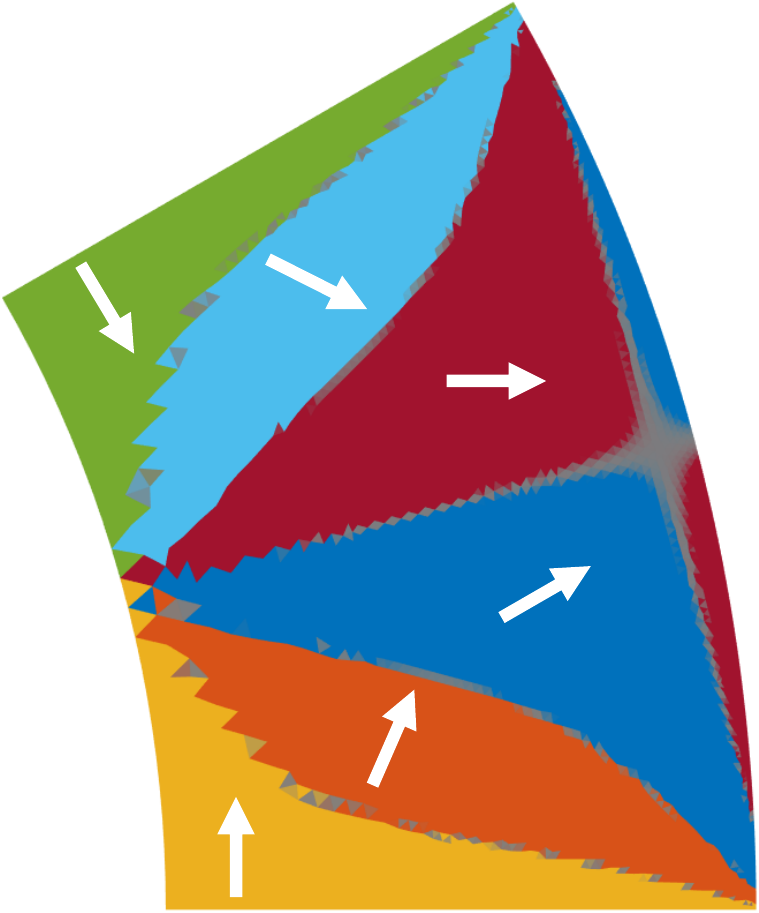} \\
		$\phi^+ = - \SI{23.1}{mWb / m} $ & $\phi^+ = \SI{7.3}{mWb / m} $ & $\phi^+ =  \SI{23.1}{mWb / m}  $ \\
  $\phi^- = - \SI{23.1}{mWb / m} $ & $\phi^- = -\SI{7.3}{mWb / m} $ & $\phi^- = \SI{23.1}{mWb / m}  $ \\
	\end{tabular}
	\caption{Non-hybridized optimized designs obtained with integers $\gamma$ and the recursive interpolation domain (results obtained with the other domains are similar). $\phi$ is normalized on \SI{1}{m} axial length.} \label{fig:Cherriere:GammaInteger}
\end{figure}

When $\gamma=0$, the flux must be maximized outwards with a positive current feeding and inwards with a negative current feeding. Then, the optimized result contains no magnets (that cannot reverse the direction of their flux) but rather conductors and steel to carry the magnetic flux.
These single-excited optimized rotors obtained with integer $\gamma$ are shown in Figure \ref{fig:Cherriere:GammaInteger}, and the color scale to interpret the materials is given in the previous Figure~\ref{fig:Cherriere:domainApplication}.

\newpage

\begin{figure}
	\centering
	\includegraphics[height=8cm]{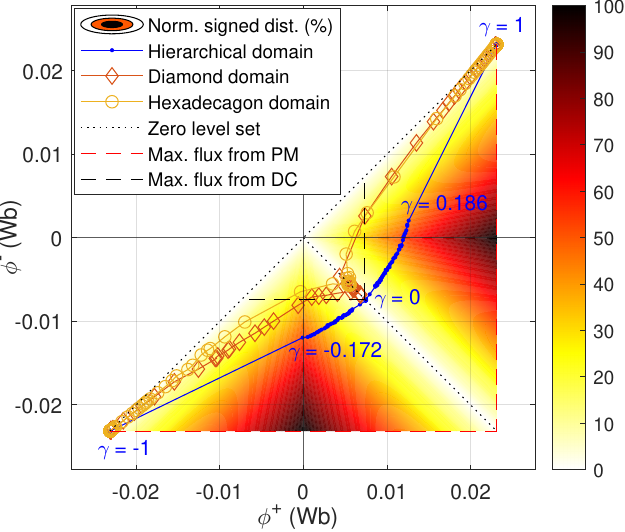}
	\caption{Pareto front obtained for $\gamma \in [-1,1]$. The step size between $\gamma$ values is $10^{-3}$.} \label{fig:Cherriere:ResultPareto}
\end{figure}

Hybrid excited rotors are obtained with intermediate (non-integer) values of $\gamma$. In particular, a "hybridization indicator", denoted $sd_0$, can be defined by the normalized signed distance to points where $|\phi^+| = |\phi^-|$. It is maximum when $|\phi^+|$ is maximum and $|\phi^-|=0$, or vice versa, corresponding to an ideal hybrid-excited rotor.  This indicator is plotted in the $(\phi^+, \phi^-)$ plane in Figure~\ref{fig:Cherriere:ResultPareto}, where the Pareto front based on recursive interpolation reaches higher values than the others, and stays between the physical lower and upper bounds of the flux.

\begin{figure}
	\centering
	\begin{tabular}{ccc}
		Hexadecagon & Diamond & Recursive domain \\ \hline \\
        $ \gamma = -0.032$ & $\gamma  = -0.068$ & $\gamma  = 0.174$ \\
		\includegraphics[width=0.25\textwidth]{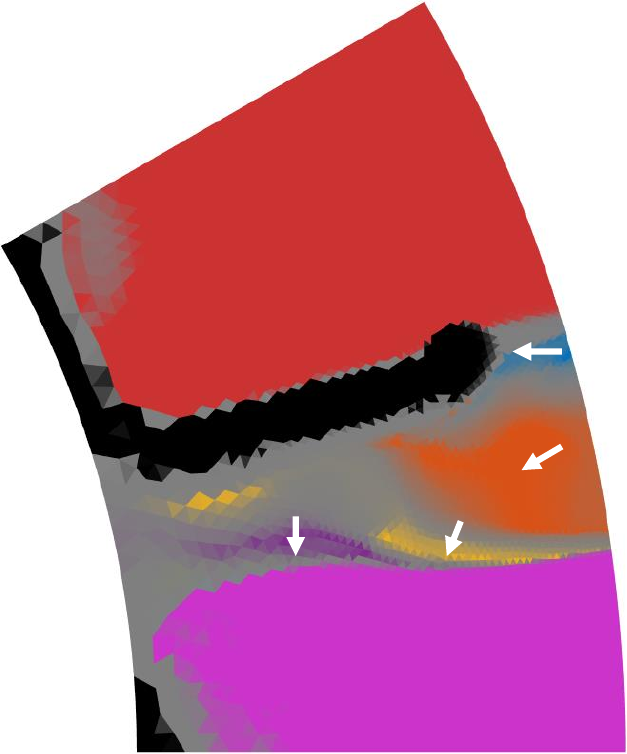} \qquad 
		& \qquad \includegraphics[width=0.25\textwidth]{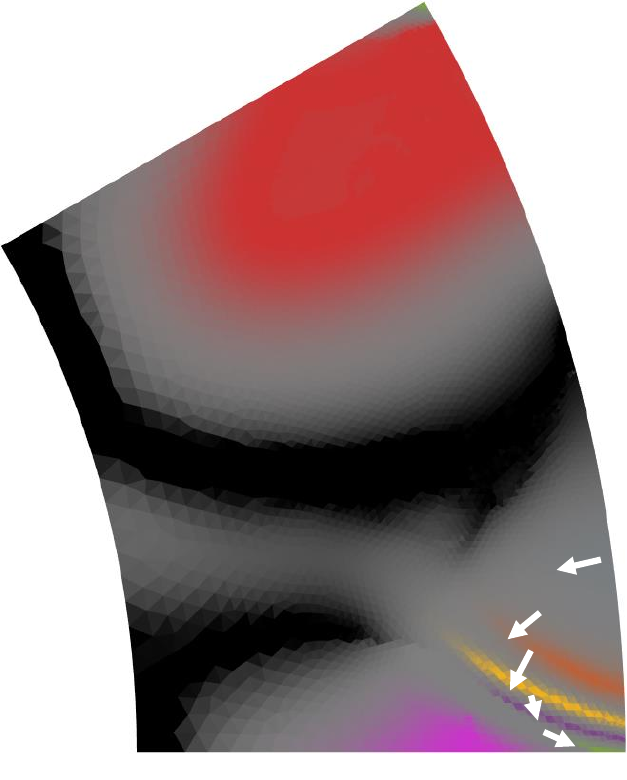}  \qquad
		& \qquad \includegraphics[width=0.25\textwidth]{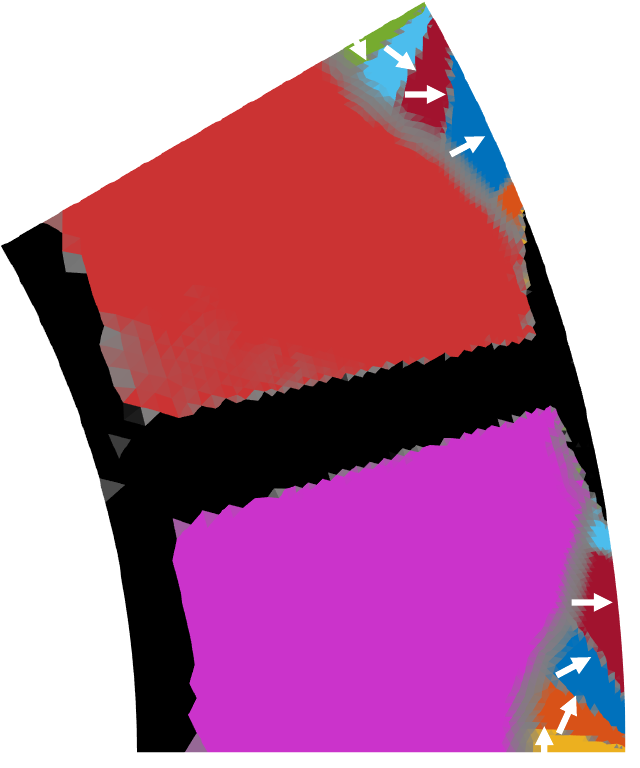} \\
		$\phi^+ = -\SI{0.3}{mWb/m}$ & $\phi^+ = -\SI{0.6}{mWb/m}$ & $\phi^+ = \SI{12.1}{mWb/m}$ \\
    $\phi^- = \SI{-6.5}{mWb/m} $ & $\phi^- = -\SI{8.1}{mWb/m} $ & $\phi^- = \SI{0.4}{mWb/m} $ \\
        $sd_0 = \SI{26}{\%}$ & $sd_0 = \SI{32}{\%}$ & $sd_0 = \SI{50}{\%}$ \\
	\end{tabular}
	\caption{Optimized designs that maximize the hybridization indicator $sd_0$ for the three different interpolation domains.}   \label{fig:Cherriere:GammaBest}
\end{figure}

 The best structures according to $sd_0$ indicator are plotted in Figure~\ref{fig:Cherriere:GammaBest}. We note that the recursive interpolation domain leads to a more performing hybrid-excited rotor. Moreover, this design is symmetric and contains fewer "gray" intermediate materials, highlighting the non-convergence in the optimizations relying on the other domains. The structure obtained by recursive interpolation is close to a machine with parallel hybrid excitation \cite{Cherriere:Hwang2018}, in which the flux created by the windings is in parallel with the flux created by the magnets. This structure is more controllable than serial double-excitation, in which the flux generated by the windings passes through the magnet, whose permeability is very low. Once again, the remaining magnets are arranged in a Halbach array, which maximizes their flux.

\section{Conclusion}

The recursive interpolation proposed in this article offers more flexibility than the interpolation domains currently available in the literature and can be adapted to a wide variety of applications. In particular, the formalism is suitable for the topology optimization of double-excitation machines, which is original. 

Further work will focus on the rotor's mechanical integrity and include additional alternating current conductors for applications such as stators of hybrid-excited flux-switching machines \cite{Cherriere:Gaussens2014}. The long-term objective is to build a comprehensive toolbox for optimizing electromagnetic actuators without initial information on the geometry and requiring as few tuning parameters as possible.

%%%%%%%%%%%%%%%%%%%%%%%%%%%%%%%%%%%%%%%%%%%%%%%%%%%%%%%%%%%%%%%%%%%%%%

%%%%%%%%%%%%%%%%%%%%%%%%%%%%%%%%%%%%%%%%%%%%%%%%%%%%%%%%%%%%%%%%%%%%%%

\end{document}